\documentclass{amsart}
\usepackage{amsfonts}
\usepackage{color}
\usepackage{amsthm}
\usepackage{amscd}

\input{tcilatex}
\input tcilatex

\begin{document}
\author{Joel E. Cohen$^{1}$ and Thierry E. Huillet$^{2}$}
\address{$^{1}$Laboratory of Populations, Rockefeller University and Columbia
University, 1230 York Avenue, Box 20, New York, NY 10065, USA; Earth
Institute and Department of Statistics, Columbia University, New York 10027,
USA; Department of Statistics, University of Chicago, Chicago IL 60637, USA;
Email: cohen@rockefeller.edu; $^{2}$Laboratoire de Physique Th\'{e}orique et
Mod\'{e}lisation\\
CY Cergy Paris University, CNRS UMR-8089\\
Site de Saint Martin\\
2 avenue Adolphe-Chauvin\\
95302 Cergy-Pontoise, France\\
E-mail: thierry.huillet@cyu.fr}
\title[Taylor's law for population models]{Taylor's law for some infinitely divisible probability distributions from population models}

\begin{abstract}
In a family of random variables, Taylor's law or Taylor's power law of
fluctuation scaling is a variance function that gives the variance $\sigma
^{2}>0$ of a random variable (rv) $X$ with expectation $\mu >0$ as a power
of $\mu$: $\sigma ^{2}=A\mu ^{b}$ for finite real $A>0,\ b$ that are the
same for all rvs in the family. Equivalently, TL holds when $\log \sigma
^{2}=a+b\log \mu ,\ a=\log A$, for all rvs in some set. Here we analyze the
possible values of the TL exponent $b$ in five families of infinitely
divisible two-parameter distributions and show how the values of $b$ depend
on the parameters of these distributions. The five families are
Tweedie-Bar-Lev-Enis, negative binomial, compound Poisson-geometric,
compound geometric-Poisson (or P\'{o}lya-Aeppli), and gamma distributions.
These families arise frequently in empirical data and population models, and 
they are limit laws of Markov processes that we exhibit in each case.

\end{abstract}
\maketitle

\section*{Keywords}

\noindent branching process; compound distribution; infinite divisibility;
Markov process; Ornstein-Uhlenbeck process; P\'{o}lya-Aeppli distribution;
power law; self-decomposability; Taylor's law; Taylor's power law; Tweedie
family; variance function

\section*{Declarations of interest}
The authors have no conflicts of interest associated with this paper.

\section*{Data availability statement}
There are no data associated with this paper.

\section{Introduction}

For practical purposes like pest control, conservation, and evaluating the
yields of alternative agricultural treatments, agronomists, ecologists, and
agricultural statisticians invented a concept today called the ``variance
function.'' In a set of samples of some quantitative observations, the
sample variance function describes how the sample variance varies from
sample to sample as a function of the mean of each sample \cite
{Bartlett1936, Beall1939, Beall1942, Bliss}.

Over several decades in the mid-twentieth century, multiple ecologists,
apparently independently, observed that the log of the sample variance was
approximated well by a linear function of the log of the sample mean, and
that this model of the sample variance function was superior to several
alternatives \cite{Bliss,Fracker,Hayman,Taylor}. The sample variance
function that specifies the log sample variance as a linear function of the
log sample mean has become known as Taylor's law (TL), after the last
ecologist who discovered it \cite{Taylor}. The mathematically equivalent
statement that the sample variance is approximately proportional to some
power of the sample mean has become known as Taylor's power law, or
fluctuation scaling in physical applications. Thousands of empirical
examples of TL have been found in many different fields of science and
finance \cite{Eisler, Taylor2019}.

In theoretical studies of TL, a set of samples is modeled by a family (set)
of random variables (rvs) $\{X(\theta )\}_{\theta \in \Theta }$ indexed by
some scalar- or vector-valued label $\theta \in \Theta $ in an index set $%
\Theta $ {with at least two elements}. The sample mean of empirical studies
is replaced by the population mean $0<EX(\theta ):=\mu (\theta )\leq \infty $%
, and the sample variance of empirical studies is replaced by the population
variance $0<$ \text{Var}$X(\theta ):=\sigma ^{2}(\theta )\leq \infty $.
Then, for finite real $a,b$ {which may or may not depend on $\theta $,} 
and for all $\theta \in \Theta $, TL asserts that $\log \sigma ^{2}(\theta
)=a+b\log \mu (\theta )$ (log-linear form) or $\sigma ^{2}(\theta )=A[\mu
(\theta )]^{b},\ A=\exp (a)$ (power form) or 
\begin{equation}
b=\frac{\log \sigma ^{2}(\theta )-a}{\log \mu (\theta )},\quad \forall
\theta \in \Theta .  \label{eq:TL}
\end{equation}

The connections between the TL parameters $a,b$ and the distributions of $%
X(\theta )$ with moments $\mu (\theta ),\sigma ^{2}(\theta )$ have long been
of interest, and have sometimes been misunderstood. The best known example
is the single-parameter family of Poisson$(\theta )$ distributions with
expectation $\theta \in (0,\infty )$, which are widely used as a model of
pure randomness in integer counts. Here $\sigma ^{2}(\theta )=\mu (\theta
)=\theta $, so TL holds with $a=0,\ b=1$ {\ for every $\theta \in (0,\infty
) $.} It does not follow that if TL holds with $a=0,\ b=1$, then any of the
distributions is Poisson$(\theta )$ \cite{CohenTE}. If $a\neq 0$ or $b\neq 1$%
, then it is safe to conclude that at least one of the distributions is not
a Poisson$(\theta )$ distribution.

Another well known single-parameter family of distributions is the negative
exponential family Exp$(\theta )$, $\theta \in (0,\infty )$. Here $\sigma
^{2}(\theta )=[\mu (\theta )]^{2}$, so TL holds with $a=0,\ b=2${,
regardless of $\theta $}. Again, it does not follow that if TL holds with $%
a=0,\ b=2$, then any of the distributions is Exp$(\theta )$.

Here we analyze the possible values of the TL exponent $b$ in five families
of infinitely divisible two-parameter distributions and show how the values
of $b$ depend on the parameters of these distributions. The shape of $%
b\left( \theta \right) $ appears to be specific to each family. A rv $X$ is
defined to be infinitely divisible if and only if, for every positive
integer $n$, there exist $n$ independent and identically distributed (iid)
rvs $X_{n,1},\ldots ,X_{n,n}$ such that $X_{n,1}+\cdots +X_{n,n}\overset{d}{=%
}X$, where $\overset{d}{=}$ means ``has the same distribution as''. 
If $X$ is infinitely divisible and 
$E(e^{itX_{n,1}})=E(e^{itX})^{1/n}$
is the characteristic function
of (say) $X_{n,1},$ then for all $t\in R$\,
\begin{equation*}
\exp \left\{ -n\left( 1-E(e^{itX_{n,1}})\right) \right\} \rightarrow E\left(
e^{itX}\right) \text{ as }n\rightarrow \infty \text{,}
\end{equation*}
showing that $X$ is a weak limit of a compound Poisson sequence. 
More general iid sequences (of size $k_{n}\rightarrow \infty$)
converging weakly to $X$ can be found in Theorem 5.2 of \cite{Steutel}.

If $X$ is infinitely divisible and satisfies Taylor's power law with
exponent $b$, then each of its constitutive summands $X_{n,m}$ also
satisfies Taylor's power law with the same exponent $b$. To see this,
suppose every rv in a family $\{X(\theta )\}_{\theta \in \Theta }$, $\Theta
\neq \emptyset $, satisfies 
\begin{equation}
X(\theta )\overset{d}{=}X_{n,1}(\theta )+\cdots +X_{n,n}(\theta ),\quad n\ge
1,  \label{eq:sumiid}
\end{equation}
where $X_{n,1}(\theta ),\ldots ,X_{n,n}(\theta )$ are iid as 
$X_{n,1}(\theta )$. Then $X(\theta )$ satisfies TL \eqref{eq:TL} with
exponent $b$ if and only if $X_{n,1}(\theta )$ satisfies TL \eqref{eq:TL}
with the same exponent $b$. Indeed, taking the expectation of both sides of %
\eqref{eq:sumiid} gives $\mu _{X}(\theta )=n\mu _{X_{n,1}}(\theta )$ and $%
\sigma _{X}^{2}(\theta )=n\sigma _{X_{n,1}}^{2}(\theta )$. If $\sigma
_{X}^{2}(\theta )=A[\mu _{X}(\theta )]^{b}$, then $n\sigma
_{X_{n,1}}^{2}(\theta )=An^{b}(\mu _{X_{n,1}}(\theta ))^{b}$ or $\sigma
_{X_{n,1}}^{2}(\theta )=An^{b-1}(\mu _{X_{n,1}}(\theta ))^{b}$, which is TL
with the same exponent $b$. The coefficient $A$ of $X(\theta )$ leads to the
coefficient $A_{n}=An^{b-1}$ of $X_{n,1}(\theta )$. The converse is obvious.

The five families to be analyzed here are Tweedie{-Bar-Lev-Enis}, negative
binomial, compound Poisson-geometric, compound geometric-Poisson (or
P\'{o}lya-Aeppli), and gamma distributions. (Kendal \cite{Kendal} gives another
interesting family.) These families arise frequently in empirical data and
population models. They all have the special property that
it is possible to express the dependence of the mean and variance 
on their parameters in such a way that
$a=0$ or $\sigma
^{2}(\theta )=[\mu (\theta )]^{b},\ \forall \theta \in \Theta $ or 
\begin{equation}
b=\frac{\log \sigma ^{2}(\theta )}{\log \mu (\theta )},\quad \forall \theta
\in \Theta .  \label{eq:TH}
\end{equation}
Non-zero $a$ can arise from the rescaling described in section \ref
{sec:scaling}.

In addition to being infinitely divisible, some of these distributions are
also self-decomposable (SD). A distribution or rv $X$ is defined to be SD
if, for every $c\in (0,1)$, there exists an independent rv $X_{c}$ such that 
\begin{equation}  \label{eq:SD}
X(\theta) \overset{d}{=} c(\theta)X(\theta) + X_{c(\theta)}(\theta).
\end{equation}
If $X$ takes values in $\mathbb{N}_{0}:=\{0,1,2,\ldots \},$ $cX$ is to be
interpreted as the $c$-Bernoulli thinning of $X$, which is defined as the
sum of $X$ iid Bernoulli rvs with success probability $c$. 

If every rv in a family $\{X(\theta )\}_{\theta \in \Theta }$, $\Theta \neq
\emptyset $, satisfies TL \eqref{eq:TL} and \eqref{eq:SD}, then the family $%
\{X_{c}(\theta )\mid \theta \in \Theta \}$ also satisfies TL with the same
value of $b$ as the family $\{X(\theta )\mid \theta \in \Theta \}$; and
conversely. To see this, take the expectation of both sides of \eqref{eq:SD}%
. Then $\mu (\theta )=c(\theta )\mu (\theta )+EX_{c(\theta )}(\theta )$,
hence $\mu (\theta )=EX_{c(\theta )}{(\theta )}/(1-c(\theta ))$. Taking the
variance of both sides of \eqref{eq:SD} gives $\sigma ^{2}(\theta )=c(\theta
)^{2}\cdot \sigma ^{2}(\theta )+\text{Var}X_{c(\theta )}(\theta )$, hence,
using TL for the second equality, 
\begin{align*}
\text{Var}X_{c(\theta )}(\theta )& =\sigma ^{2}(\theta )\cdot (1-c(\theta
)^{2}) \\
& =A[\mu (\theta )]^{b}\cdot (1-c(\theta )^{2}) \\
& =A[EX_{c(\theta )}{(\theta )}/(1-c(\theta ))]^{b}\cdot (1-c(\theta )^{2})
\\
& =\{A(1-c(\theta )^{2})(1-c(\theta ))^{-b}\}[EX_{c(\theta )}{(\theta )}%
]^{b}.
\end{align*}
This is TL with the same $b$, as claimed. The coefficient $A$ of $X(\theta )$
leads to the coefficient $A(1-c(\theta )^{2})(1-c(\theta ))^{-b}$ of $%
X_{c(\theta )}(\theta )$. The proof of the converse is routine, given the
above equalities.

Every SD rv is unimodal. If $X$ is SD and not discrete-valued, then it has a
density. If $X$ is SD and continuous, it is a weak limit of a
continuous-time L\'{e}vy-driven Ornstein-Uhlenbeck (OU) process. If $X$ is
SD and $\mathbb{N}_{0}$-valued, it is a weak limit of a pure-death branching
process with immigration. See the Appendix.

Thus we consider three kinds of infinitely divisible distributions with two
parameters. First, ``bare'' infinitely divisible distributions arise as
specific limit laws, as we will describe. 
Second, infinitely divisible distributions that are compound Poisson
geometric are the limit laws of immigration processes with total disasters. Third, infinitely divisible distributions that are
SD are the limit laws of a pure-death branching process with immigration (if
discrete) or of a L\'{e}vy-driven Ornstein-Uhlenbeck process (if
continuous). In the latter two cases, there is a balance between events of birth 
(or immigration)  and events of
death, resulting in an equilibrium distribution. {For some of the
families under study, depending on the parameter range, the rv can switch
from ``bare'' infinitely divisible to SD, perhaps suggesting some kind of
phase transition. } {There may be other examples of infinitely divisible
distributions that fall into these three categories, and it is of interest
to know whether they would satisfy Taylor's law. }

\section{Tweedie{-Bar-Lev-Enis} distributions}

The definition of Tweedie{-Bar-Lev-Enis} distributions is elaborate \cite
{JB, JBB, BLE1986, BLS1987} and will not be attempted here. An accessible
expository account, with historical background, is \cite{BL2020}. For
brevity, we will refer to these distributions and associated random
variables as the TweBLE family.

With parameter $\alpha \in [ -\infty ,0) \cup \{ 0\} \cup (0,1) \cup ( 1,2) \cup \{2\} $ 
and 
\begin{equation}
k\left( \theta \right) =\frac{1-\alpha }{\alpha }\left( \frac{\theta }{
1-\alpha }\right) ^{\alpha }  \label{ktheta}
\end{equation}
for the values of parameter $\theta $ for which 
$\left( \frac{\theta }{1-\alpha }\right) ^{\alpha }$ 
is well-defined, the bilateral probability
Laplace-Stieltjes transform (PLSt) of a TweBLE random variable (rv) \cite
{JB, JBB} is 
\begin{equation*}
\Phi ( \lambda ) =e^{-[ k( \theta +\lambda )-k( \theta ) ] }.
\end{equation*}
The expression of $k( \theta ) $ is extended to cover the
boundary cases $\alpha =-\infty $ and $\alpha =0$, respectively, by 
$1-e^{-\theta }$ and $\log \theta $. TweBLE rvs include Poisson ($\alpha
=-\infty $), compound Poisson-gamma ($\alpha \in ( -\infty ,0) $)
(which is the sum of a Poisson-distributed number of iid gamma rvs all
independent of the Poisson), negative exponential ($\alpha =0$), tempered
one-sided stable ($\alpha \in ( 0,1) $), tempered two-sided
extreme stable ($\alpha \in ( 1,2) $) and Gaussian ($\alpha =2$)
rvs.

Considering (\ref{ktheta}), the log-Laplace transform (LLt) is 
\begin{equation*}
L( \lambda ) :=-\log \Phi ( \lambda ) =\frac{1-\alpha 
}{\alpha }\left[ \left( \frac{\theta +\lambda }{1-\alpha }\right) ^{\alpha
}-\left( \frac{\theta }{1-\alpha }\right) ^{\alpha }\right] .
\end{equation*}
Hence the mean and variance of a TweBLE rv are 
\begin{eqnarray*}
\mu &=&L^{\prime }( 0) =\left( \frac{\theta }{1-\alpha }\right)
^{\alpha -1}, \\
\sigma ^{2} &=&-L^{^{\prime \prime }}( 0) =\left( \frac{\theta }{%
1-\alpha }\right) ^{\alpha -2}.
\end{eqnarray*}
All TweBLE distributions satisfy $\sigma ^{2}=\mu ^{b}$, which is TL with $%
a=0$. Taylor's law holds with exponent 
\begin{equation}
b=\frac{\log \sigma ^{2}}{\log \mu }=\left( 2-\alpha \right) /\left(
1-\alpha \right) .  \label{bofalpha}
\end{equation}
The extensions of TweBLE rvs to $\alpha =-\infty $ and $\alpha =0$ must have 
$b=1$ and $b=2$, respectively, to be consistent, and these values of $b$
follow directly from the moments of Poisson and exponential distributions,
respectively. Although a TweBLE rv has two parameters $\alpha ,\ \theta $,
the TL exponent $b$ depends only on $\alpha $, unlike the two-parameter
distributions in the following sections.

The graph of $b$ versus $\alpha $ when $\alpha \in \left[ 0,2\right] $
diverges at $\alpha _{c}:=1$ and has two hyperbolic branches:

- one increasing from $\alpha =-\infty $ to $\alpha =1^{-}$, with $b$
ranging from $1^{+}$ (Poisson) to $b=+\infty $, through $b=2$ (gamma) when $%
\alpha =0.$

- one increasing from $\alpha =1^{+}$ to $\alpha =2$, with $b$ ranging from $%
b=-\infty $ to $b=0$ ({Gaussian} when $b=0,\ \alpha =2$).

A TweBLE rv has a negative TL exponent $b<0$ only when its support is the
whole real line, unlike some of the following examples, where $b<0$ can
occur when the support is the nonnegative half line or the nonnegative
integers.

No TweBLE distribution has $b\in \left( 0,1\right) $ because there is no
non-degenerate TweBLE PLSt with $\theta >2$. Hence a population model that
obeys TL with $b\in \left( 0,1\right) $ cannot have the distribution of, and
cannot be explained by using, a TweBLE rv. Examples of TL with $b\in \left(
0,1\right) $ include a multiplicative population process in a Markovian
environment \cite[p. 33]{Cohen2014} and an infinity of examples with
arbitrary distributions with finite means and finite variances 
\cite[Example 3]{CohenTE}. Therefore TweBLE distributions are not sufficient
for modeling and understanding all instances of TL.

When $\alpha<1$, TweBLE rvs are infinitely divisible because 
$\Phi (\lambda ) =\exp -L( \lambda )$ where 
$L^{\prime}( \lambda ) $ is completely monotone, obeying 
$(-1) ^{n}L^{( n+1) }( \lambda ) \geq 0$
for all $\lambda >0$. 
When $\alpha \in [1,2]$, TweBLE rvs are tempered $\alpha$-stable.
In all cases, the L\'{e}vy jump measure can be found in 
\cite{BLBL1992}.
Moreover,
\cite[Example 3]{BLBL1992} shows that they are SD as well except when $%
\alpha \in ( -\infty ,0) $. For TweBLE rvs that are SD, the PLSt
is 
\begin{equation}
\Phi \left( \lambda \right) =\exp \left\{ \int_{0}^{\lambda }\frac{\log \Phi
_{0}\left( \lambda ^{\prime }\right) }{\lambda ^{\prime }}d\lambda ^{\prime
}\right\}  \label{Phi}
\end{equation}
for some PLSt $\Phi _{0}( \lambda ) $ of an infinitely divisible
rv \cite[Theorem 2.9, Eq. 2.12]{Steutel}.

TweBLE distributions as defined in this Section belongs to the class
of exponential families \cite{BL2020}. 
All such TweBLE distributions
satisfy $\sigma ^{2}(\theta )=[\mu (\theta )]^{b}$, which is TL
with $A=1,\ a=0$, and $b$
given by \eqref{bofalpha}. Introducing the scaling %
\eqref{eq:rescale} below, we switch to the class of exponential dispersion
models, as defined by Jorgensen \cite{JBB}. For this new class, $\sigma
^{2}(\theta )=A[\mu (\theta )]^{b}$ with $A=e^a$ 
and $a\neq 0$, independent of $b$.

\section{Negative binomial distributions}

We shall show that the next two-parameter 
probability distribution families can be formulated in such a way that they
also obey TL with the TweBLE family's special property that
$a=0$, and we shall express $b$ in terms of their parameters. As for
TweBLE models, the scaling \eqref{eq:rescale} will yield the 
version of TL with $a\neq 0$, independent of $b$.

With the parameters $\alpha >0$, $p\in ( 0,1),\ q:=1-p$, a negative binomial rv $X%
\overset{d}{=} NB(\alpha ,p)$ has probability mass function 
\begin{equation*}
P[X=k]=\binom{k+\alpha -1}{k}q^{\alpha }p^{k},\quad k=0,1,2,\ldots .
\end{equation*}
It has PLSt and LLt 
\begin{align*}
\Phi \left( \lambda \right) & =\left( \frac{q}{1-pe^{-\lambda }}\right)
^{\alpha }, \\
L\left( \lambda \right) & =-\log \Phi \left( \lambda \right) =-\alpha \log
q+\alpha \log \left( 1-pe^{-\lambda }\right) , \\
L^{\prime }\left( \lambda \right) & =\frac{\alpha pe^{-\lambda }}{%
1-pe^{-\lambda }},
\end{align*}
and mean and variance 
\begin{eqnarray}
\mu &=&L^{\prime }\left( 0\right) =\frac{\alpha p}{q}, \nonumber \\
\sigma ^{2} &=&-L^{^{\prime \prime }}\left( 0\right) =\frac{\alpha p}{q^{2}}%
=\mu +\frac{1}{\alpha }\mu ^{2}.\label{eq:NBvariance}
\end{eqnarray}
Overdispersion holds, in the sense that, unlike the Poisson distribution,
where the variance equals the mean, the negative binomial variance strictly
exceeds the mean. The smaller $\alpha $ is, the larger the variance compared
to the mean. In ecology, overdispersion in population counts is sometimes
interpreted as a resulting from heterogeneous field conditions
(environmental variation making some places more favorable for individuals
of the population than other places) or from aggregation or clustering
(individuals behaving in a way that brings them near other individuals) \cite
{Bliss,Taylor}.

Taylor's law holds with exponent
\begin{equation}
b=\frac{\log \sigma ^{2}}{\log \mu }=1-\frac{\log q}{\log \left( \frac{%
\alpha p}{q}\right) }.  \label{eq:NBDb}
\end{equation}
Then 
\begin{align*}
b&<1\Leftrightarrow \frac{\alpha p}{q}<1\Leftrightarrow p<p_{c}:=1/\left(
1+\alpha \right) ,\\
b &<0\Leftrightarrow \frac{\log q}{\log \left( \frac{\alpha p}{q}\right) }%
>1\Leftrightarrow q^{2}+\alpha q-\alpha <0 \\
&\Leftrightarrow q<q_{0}:=\frac{-\alpha +\sqrt{\alpha ^{2}+4\alpha }}{2}%
\Leftrightarrow p>p_{0}:=\frac{2+\alpha -\sqrt{\alpha ^{2}+4\alpha }}{2}>0.
\end{align*}
We always have 
\begin{equation*}
p_{0}<p_{c}=\frac{1}{1+\alpha }\text{.}
\end{equation*}
So {for any fixed $\alpha $,} the graph of $b$ against $p$ has two branches:

- one is concave and decreasing from $b=1$ to $b=-\infty $ as $p$ varies
from $0^{+}$ to $p_{c}^{-}$ as a diverging point according as $b\rightarrow
-\infty $ as $p\rightarrow p_{c}^{-}$, passing through $0$ as $p=p_{0}.$

- one is varying from $b=+\infty $ as $p\rightarrow p_{c}^{+}$ to $b=2$ as $%
p\rightarrow 1^{-}$ possibly passing through a minimum $b_{\min }>1$ with $%
b_{\min }<2$ ($=2$) if and only if $\alpha >1$ (respectively $\alpha \leq 1$%
). Hence 
\begin{eqnarray*}
b &\in &\left( 0,1\right) \text{ if }0<p<p_{0}, \\
b &\leq &0\text{ if }p_{0}\leq p<p_{c}, \\
b &>&b_{\min }\text{ if }p>p_{c}.
\end{eqnarray*}
The range $\left( 1,b_{\min }\right) $ for $b$ is excluded.

If $\alpha =1$, then $b_{\min}=2$ and $b$ cannot fall in the interval $%
\left(1,2\right)$. In this case, the negative binomial distribution reduces
to the geometric distribution, and $p_{c}=1/2$ and $p_{0}=\left( 3-\sqrt{5 }%
\right) /2=\phi ^{2}$, where $\phi =(\sqrt{5}-1)/2$ is a golden ratio.

The formula \eqref{eq:NBDb} identifies a manifold of $(\alpha,\ p)$ points on
which the TL exponent $b$ is constant. To maintain a constant $b$, %
\eqref{eq:NBDb} requires both $\alpha$ and $p$ to vary. This result is
surprising in light of the above negative binomial variance function
\eqref{eq:NBvariance}, $%
\sigma ^{2} =\mu +\alpha^{-1}\mu ^{2}$. For suppose $\alpha$ is constant.
Then as $\mu\to 0$, $\mu^2\ll\mu$ and asymptotically $\sigma ^{2}$ is
proportional to $\mu$, which is asymptotically TL with $a=0,\ b=1$. However,
as $\mu\to \infty$, $\mu^2\gg\mu$ and asymptotically $\sigma ^{2}$ is
proportional to $\alpha^{-1}\mu ^{2}$, which is asymptotically TL with $%
e^a=\alpha^{-1},\,\, b=2$. Moreover, $\log \sigma^{2}$ is a strictly convex
increasing function of $\log \mu$ \cite[Supp.Mat., lines 29-42]{CPL},
instead of a linear function as in TL. If the expectation $\mu$ varies over
only a small range in the data or the family of rvs, then $\log\sigma^{2}$
can be closely approximated by (may be statistically indistinguishable from,
in the case of empirical data) a locally linear function of $\log\mu$ with a
local slope that increases smoothly from 1 (for very small $\mu$) to 2 (for
very large $\mu$). A wide range of means $\mu$ would be required to display
the curvature in the relation of $\log \sigma^{2}$ to $\log \mu$ when the
underlying distributions are negative binomial with constant $\alpha$.

These observations illustrate that the form of the variance function
of a two-parameter family of probability distributions depends
strongly on how the parameters and their relation are constrained.
We illustrate this point again in analyzing the gamma distribution below.

The negative binomial distribution is discrete and infinitely divisible
(compound Poisson) because its probability generating function (pgf) $\phi
\left( z\right) :=E\left( z^{X}\right) $ is 
\begin{equation}
\phi \left( z\right) =\left( \frac{q}{1-pz}\right) ^{\alpha }=\exp \{-\left[
-\alpha \log q\right] \left( 1-c\left( z\right) \right) \},
\label{eq:negbinpgf}
\end{equation}
{where} $c( z) =\log ( 1-pz) /\log q$ is the
pgf of a logarithmic series distribution for its clusters' sizes.

The negative binomial distribution is also SD because its pgf may be written
as (compare with \eqref{Phi} in the continuum) 
\begin{equation*}
\phi \left( z\right) =\exp \left\{ -r\int_{z}^{1}\frac{1-h\left( z^{\prime
}\right) }{1-z^{\prime }}dz^{\prime }\right\}
\end{equation*}
for some rate $r>0$ and pgf $h( z) $ such that $h( 0) =0$ 
\cite[Theorem 4.11, Eq. 4.13]{Steutel}. Indeed, $r=\alpha p>0$ and $h( z)
=qz/(1-pz)$ is a geometric pgf. See the Appendix.

When $\alpha \in (0,1)$, then $\pi :=1-q^{\alpha }\in (0,1)$ and the pgf $%
\phi (z)$ in \eqref{eq:negbinpgf} equals 
\begin{equation*}
\phi (z)=\frac{1-\pi }{1-\pi \varphi (z)},
\end{equation*}
where 
\begin{equation*}
\varphi (z):=\frac{1-(1-pz)^{\alpha }}{1-q^{\alpha }}
\end{equation*}
satisfies $\varphi (0)=0$. This $\varphi (z)$ is a pgf because it is 
an absolutely monotone function:
it obeys $\varphi ^{(n)} ( z) \geq 0$ for all $0<z<1$
and has positive coefficients in its power series expansion. 
Indeed, the coefficient of $z^{n}$ in $\varphi (z)$ is 
\begin{eqnarray*}
[ z^{n}]\varphi ( z) &=&\frac{p^{n}}{1-q^{\alpha }}(
-1) ^{n-1}\frac{( \alpha ) _{n}}{n!} \\
&=&\frac{p^{n}}{1-q^{\alpha }}\frac{\alpha [ \overline{\alpha }]
_{n-1}}{n!}, \ n\geq 1
\end{eqnarray*}
where $\left( \alpha \right) _{n}:=\Gamma \left( \alpha +1\right) /\Gamma
\left( \alpha +1-n\right) =\alpha \left( \alpha -1\right) ...\left( \alpha
-n+1\right) $ is the falling factorial and $\left[ 
\overline{\alpha }\right] _{n}:=\overline{\alpha }\left( \overline{\alpha }%
+1\right) ...\left( \overline{\alpha }+n-1\right) ,$ $n\geq 1$, is the
rising factorial of $\overline{\alpha }$ with $\left[ \overline{\alpha }%
\right] _{0}:=1$, $\overline{\alpha }=1-\alpha $.

Because the negative binomial distribution is SD for all $\alpha >0$, it
follows from our Appendix that the negative binomial distribution is the
limit law of a pure-death branching process with immigration. But when $%
\alpha \in ( 0,1)$, the negative binomial distribution is also
compound-geometric, so it is also the limit law of a Markov chain with total
disasters, as defined in the next section.

\section{Compound Poisson-geometric distributions}

With parameters $\alpha >0,\ p\in ( 0,1) $, $q:=1-p,$ a compound
Poisson-geometric rv is the sum of $N$ iid Poisson($\alpha$)-distributed rvs,
where $N+1$ is geometrically distributed with parameter $p$. 
The PLSt and LLt of a compound Poisson-geometric rv are 
\begin{align*}
\Phi ( \lambda ) & =\frac{q}{1-pe^{-\alpha ( 1-e^{-\lambda
}) }}, \\
L( \lambda ) & :=-\log \Phi ( \lambda ) =-\log q+\log
\left( 1-pe^{-\alpha ( 1-e^{-\lambda }) }\right) .
\end{align*}
The mean and variance are 
\begin{eqnarray*}
\mu &=&L^{\prime }( 0) =\frac{\alpha p}{q}, \\
\sigma ^{2} &=&-L^{^{\prime \prime }}( 0) =\frac{\alpha p(
\alpha +q) }{q^{2}}=\frac{\alpha +q}{q}\mu >\mu .
\end{eqnarray*}
Overdispersion holds. 
Taylor's law holds with exponent
\begin{equation}  \label{eq:PGb}
b=\frac{\log \sigma ^{2}}{\log \mu }=1+\frac{\log \left( 1+\frac{\alpha }{q}%
\right) }{\log \left( \frac{\alpha p}{q}\right) }.
\end{equation}
The TL exponent $b$ diverges when $q=q_{c}:=\alpha /( 1+\alpha ) $
or $\alpha =q/p$. It vanishes when $q=q_{0}$, defined as the solution in $%
\left( 0,1\right) $ of 
\begin{equation*}
f(q):=\left( 1+\alpha \right) q^{2}-q\alpha \left( 1-\alpha \right) -\alpha
^{2}=0.
\end{equation*}
Further, because $f(q_{c})=-\alpha ^{2}/(1+\alpha )<0$, we have $q_{0}>q_{c}$%
. So the graph of $b$ has two branches {for a given fixed $\alpha$}:

- one is increasing and convex from $b=2$ to $b=\infty $ as $q$ varies from $%
0^{+}$ to $q_{c}^{-}$.

- one is increasing from $b=-\infty $ to $b=1$ as $q$ increases from 
$q_{c}^{+}$ to 1, passing through $0$ as $q=q_{0}>q_{c}$. Hence 
\begin{eqnarray*}
b &\in &( 2,\infty ) \text{ if }0<q<q_{c}, \\
b &<&1\text{ if }q>q_{c}.
\end{eqnarray*}
The TL exponent $b$ is excluded from $( 1,2) $.

If $\alpha >0$ is small enough that $e^{-\alpha }<(1-\alpha )/\alpha$ 
(e.g., if $\alpha <0.659)$ and if $p\geq p_*:=\alpha
/(1-e^{-\alpha }\alpha )\in (0,1)$, then the compound Poisson-geometric rv
is SD.

To prove this claim, we observe that, with $\varphi ( z) =e^{-\alpha ( 1-z)}$
\begin{equation*}
\phi ( z) =\frac{q}{1-p\varphi (z)}=e^{-( -\log q) [ 1-c( \varphi ( z) ) ] }
\end{equation*}
where $c( z) $ is the pgf of a logarithmic series distribution. So $\phi (
z) $ is compound Poisson with clusters' size pgf $c( \varphi ( z) ) $.

The probability mass function $P[ X=k] =[ z^{k}] \phi ( z) $ 
of a compound Poisson-geometric rv $X$ is in principle
explicitly given by the Faa di Bruno formula for compositions of two pgfs 
\cite[p. 146]{Comtet}. It involves the ordinary Bell polynomials in the
variables $b_{k}=[z^{k}] \varphi ( z)$.

Next, we try to write the pgf $\phi (z) $ of $X$ in the form 
\begin{equation*}
\phi ( z) =\frac{q}{1-p\varphi (z)}=\exp \{-r\int_{z}^{1}\frac{1-h(
z^{\prime }) }{1-z^{\prime }}dz^{\prime }\}
\end{equation*}
for some rate $r>0$ and pgf $h(z)$ obeying $h(0)=0$. Forcing $h(0)=0$ yields 
$$r=\frac{p\varphi \varphi ^{\prime }( 0) }{1-p\varphi ( 0) }=\frac{p\alpha
e^{-\alpha }}{1-pe^{-\alpha }}.$$ Hence 
\begin{equation*}
h( z) =1-\frac{1}{r}( 1-z) \frac{p\varphi ^{\prime }( z) }{1-p\varphi ( z) }=%
\frac{1}{r}\frac{r( 1-p\varphi ( z) ) -p( 1-z) \varphi ^{\prime }( z) }{%
1-p\varphi ( z) }.
\end{equation*}
Denoting the numerator by $N(z)$, a sufficient condition for $h$ to be a pgf
is that 
\begin{equation*}
\lbrack z^{k}]N(z)\geq 0\text{ for all }k\geq 1,
\end{equation*}
where $[z^{k}]N(z)$ denotes the coefficient of $z^{k}$ in the power series
expansion of $N(z)$. But, with $b_{k}=e^{-\alpha }\alpha ^{k}/k!$, 
we have 
\begin{equation*}
N(z)=p\sum_{k\geq 1}z^{k}[(k-r)b_{k}-(k+1)b_{k+1}].
\end{equation*}
So $h$ is a pgf if 
\begin{equation*}
\frac{b_{k+1}}{b_{k}}=\frac{\alpha }{k+1}\leq \frac{p( k-b_{1}) }{k+1}\text{
for any }k\geq 1
\end{equation*}
or equivalently 
\begin{equation*}
\alpha \leq p( 1-b_{1}) =p( 1-e^{-\alpha }\alpha ) \text{ or }p\geq
p_{*}:=\alpha /( 1-e^{-\alpha }\alpha ) .
\end{equation*}
\noindent This completes the proof that, under the above condition, the
compound Poisson-geometric rv is SD.

A compound Poisson-geometric rv is the limiting rv of a Markov chain with
total disasters. To see this, let $(\beta _{n})_{n\geq 1}$ be an iid
sequence taking values in $\mathbb{N}_{0}:=\{0,1,2,...\}$. Let $\varphi
(z):=E(z^{\beta })$ be the pgf of the $\beta $s. The Markov chain $X_{n}$
evolves by: 
\begin{align*}
X_{n+1}& =X_{n}+\beta _{n+1}\text{ with probability }p, \\
X_{n+1}& =0\quad \quad \quad \quad \ \text{ with probability }q:=1-p.
\end{align*}
\noindent This simple population growth model alternates periods of births
of amplitude $\beta _{n+1}$ with one or more total disasters where
population size $X_{n}$ is instantaneously reset to $0$. We assume without
loss of generality that $X_{0}=0$ and that $q$ does not depend on $X_{n}=x$. 
{We define $p*X$\ to equal the product of $X$ times an independent Bernoulli
rv $B(p)$ with success parameter $p$. (This product is not to be confused
with Bernoulli thinning.) } Because $X:=X_{\infty }$ solves the
distributional equation $X\overset{d}{=}B(p)(X^{\prime }+\beta
)=:p*(X^{\prime }+\beta )$, where $(B(p),X^{\prime }\overset{d}{=}X,\beta )$
are mutually independent and $B(p)$ is a Bernoulli rv with success parameter 
$p$, the above Markov chain $\{X_{n}\mid
n=0,1,2,\ldots\} $ with total disasters is clearly ergodic. The limiting rv $%
X:=X_{\infty }$ exists and it has the compound geometric pgf 
\begin{equation*}
\phi (z)=q/(1-p\varphi (z))
\end{equation*}
(shifted to the left by one unit).
In our example, $\beta $ is Poisson-distributed with pgf 
$\varphi(z)=e^{-\alpha (1-z)}$.

\section{Compound geometric-Poisson distributions}

With parameters $\alpha >0,\ p\in ( 0,1),\ q:=1-p$, a compound geometric-Poisson (or
P\'{o}lya-Aeppli) rv is the sum of a Poisson($\alpha$)-distributed number of
geometric rvs with parameter $p$. The PLSt and LLt are

\begin{align*}
\Phi (\lambda )& =e^{-\alpha \left( 1-\frac{qe^{-\lambda }}{1-pe^{-\lambda }}%
\right) }=e^{-\alpha \frac{1-e^{-\lambda }}{1-pe^{-\lambda }}}, \\
L(\lambda )& :=-\log \Phi (\lambda )=\alpha \frac{1-e^{-\lambda }}{%
1-pe^{-\lambda }}.
\end{align*}
Hence the mean and variance are 
\begin{eqnarray*}
\mu &=&L^{\prime }(0)=\frac{\alpha }{q}>0, \\
\sigma ^{2} &=&-L^{^{\prime \prime }}(0)=\frac{\alpha (1+p)}{q^{2}}=\frac{1+p%
}{q}\mu >\mu .
\end{eqnarray*}
Its probability mass function is given by 
\begin{eqnarray*}
P[X=k] &=&e^{-\alpha }\sum_{l=1}^{k}\binom{k-1}{l-1}\frac{\alpha ^{l}}{l!}%
p^{k-l}q^{l}\text{, if }k\geq 1 \\
&=&e^{-\alpha }\text{ if }k=0
\end{eqnarray*}
Overdispersion holds. Taylor's law holds with exponent 
\begin{equation}
b=\frac{\log \sigma ^{2}}{\log \mu }=1+\frac{\log \left( \frac{1+p}{q}
\right) }{\log \left( \frac{\alpha }{q}\right) }.  \label{eq:GPb}
\end{equation}

If $\alpha \geq 2$, the graph of $b$ has only one decreasing branch as $q$
varies from $0$ to $1$, with $b\rightarrow 2^{-}$ as $q\rightarrow 0^{+}$
and $b=1$ at $q=1$. Hence $b\in[1,2]$.

If $1<\alpha<2$, the graph of $b$ increases from $b=2$ to a maximum 
less than 3 and then drops to 1 as $q$ increases from $0$ to $1$.

If $\alpha=1$, then $b=\log(1+p)=\log(2-q)$, so $b$ falls from 2 to 0
as $q$ increases from 0 to 1.

If $0<\alpha <1$, $b$ diverges when $q=q_{c}:=\alpha ,$ so $b$ has two
branches:

- one is increasing and convex from $b=2^{+}$ to $b=+\infty $ as $q$ increases
from $0^{+}$ to $q_{c}^{-}$.

- one is increasing and concave from $b=-\infty $ to $b=1$ as $q$ increases
from $q_{c}^{+}$ to $1^{-}$. This branch passes through $b=0$ when $%
q=q_{0}:=( \sqrt{\alpha ^{2}+8\alpha }-\alpha ) /2>q_{c}.$ Hence the full
range of $b$ is covered: 
\begin{eqnarray*}
b &\in &( 2,+\infty ) \text{ if }0<q<q_{c}, \\
b &\in &( -\infty ,2) \text{ if }q>q_{c}.
\end{eqnarray*}

The pgf $\phi (z)$ of a geometric-Poisson or P\'{o}lya-Aeppli rv is 
\begin{equation*}
\phi (z)=\exp \left\{ -\alpha \frac{1-z}{1-pz}\right\} =\exp \left\{
-r\int_{z}^{1}\frac{1-h\left( z^{\prime }\right) }{1-z^{\prime }}dz^{\prime
}\right\}
\end{equation*}
for some rate $r>0$ and pgf $h(z)$ obeying $h(0)=0$, only if $p>1/2$. The
necessary condition $p>1/2$ arises because $h(0)=0$ yields $r=\alpha q$ and 
\begin{equation*}
h(z)=z\frac{1-2p+p^{2}z}{(1-pz)^{2}},
\end{equation*}
which is absolutely monotone only if $%
p>1/2 $, and then $h(z)$ satisfies
$[z^{n}] h(z)\geq 0$ for all $n\geq 1$.

\section{Gamma distributions}

A gamma rv $X$ with shape parameter $\alpha >0$ 
and scale parameter $\beta >0 $ has probability law 
\begin{equation*}
P( X\in dx) =\frac{x^{\alpha -1}e^{-x/\beta }}{\Gamma (\alpha )\beta
^{\alpha }}dx,\ \ x\in (0,\infty ).
\end{equation*}
The PLSt and LLt are 
\begin{align*}
\Phi( \lambda ) & =\left( \frac{1}{1+\lambda \beta }\right) ^{\alpha }, \\
L( \lambda ) & =-\log \Phi ( \lambda ) =\alpha \log( 1+\lambda \beta ) .
\end{align*}
The mean and variance are 
\begin{eqnarray*}
\mu &=&L^{\prime }( 0) =\alpha \beta >0, \\
\sigma ^{2} &=&-L^{^{\prime \prime }}( 0) =\alpha \beta ^{2}=\frac{1}{\alpha 
}\mu ^{2}=\beta \mu .
\end{eqnarray*}
Overdispersion holds if and only if $\beta >1$. 

As a reviewer pointed out, a gamma rv can satisfy TL in multiple ways. 
For example, if we fix $b=2$, then $\sigma^2/\mu^2=1/\alpha=A$
and $a=-\log\alpha$.
If we fix $b=1$, then $\sigma^2/\mu=\beta=A$ and $a=\log \beta$.
If we fix $b=3/2$, then $\sigma^2/\mu^{3/2}=(\beta/\alpha)^{1/2}$
and $a=(1/2)\log(\beta/\alpha)$. 
Here we fix $a=0,\ A=1$ so that
TL holds with exponent 
\begin{equation}
b=\frac{\log \sigma ^{2}}{\log \mu }=1+\frac{\log \beta }{\log ( \alpha
\beta ) }.  \label{eq:GAMMAb}
\end{equation}

If $\beta =1$, then $b=1$ constant (equidispersion, or variance equal to the
mean).

If $\beta \neq 1$, the graph of $b$ versus $\alpha$ shows a singularity at $%
\alpha _{c}:=1/\beta ,$ with the full range of $b$ covered.

If $\beta <1$, the graph of $b$ versus $\alpha $ has two increasing
branches. We have $b\rightarrow 1^{+}$ if $\alpha \rightarrow
0^{+}, $ $b\rightarrow +\infty $ if $\alpha \rightarrow \alpha _{c}^{-}$ and 
$b\rightarrow -\infty $ if $\alpha \rightarrow \alpha _{c}^{+},$ $%
b\rightarrow 1^{-}$ if $\alpha \rightarrow \infty .$

If $\beta >1,$ the graph of $b$ is a symmetric image of the previous one
with respect to the horizontal line $b=0$, so with two decreasing branches$.$

If $\alpha =1$, giving an exponential distribution, then $b=2$ is a fixed
point.

The gamma rv is SD because \cite[p. 4, Theorem 2.9]{Steutel} 
\begin{equation*}
\Phi (\lambda )=\exp \left\{ \int_{0}^{\lambda }\frac{\log \Phi _{0}\left(
\lambda ^{\prime }\right) }{\lambda ^{\prime }}d\lambda ^{\prime }\right\}
\end{equation*}
for some PLSt $\Phi _{0}(\lambda )$ of an infinitely divisible rv, which is
here found to be 
\begin{equation*}
\Phi _{0}(\lambda )=e^{\alpha \beta \lambda /( 1+\lambda \beta )
}=e^{-\alpha ( 1-1/( 1+\lambda \beta ) ) }.
\end{equation*}
This $\Phi _{0}(\lambda )$ is the PLSt of a compound Poisson$( \alpha
) $ rv with Exp$(\beta )$ distribution for the random size $\Delta $
of its clusters. Consider the Ornstein-Uhlenbeck process 
\begin{equation*}
dX_{t}=-X_{t}dt+d\mathcal{L}_{t},\ \ X_{0}=0,
\end{equation*}
driven by the L\'{e}vy-process $\mathcal{L}_{t}$ (here a rate-$\alpha $
compound-Poisson exponential process) for which 
\begin{equation*}
Ee^{-\lambda \mathcal{L}_{t}}=\Phi _{0}(\lambda )^{t},\ \ t\geq 0.
\end{equation*}
With $\phi _{\Delta }(\lambda )={E}e^{-\lambda \Delta }=(1+\lambda
\beta )^{-1}$, 
\begin{align*}
\Phi _{t}(\lambda )&:= {E}(e^{-\lambda X_{t}}\mid X_{0}=0)=\exp
\left\{ -\alpha \int_{0}^{t}(1-\phi _{\Delta }(\lambda e^{-s}))ds\right\}
=\left( \frac{1+\lambda \beta }{1+\lambda \beta e^{-t}}\right) ^{-\alpha } \\
&\rightarrow \left( 1+\lambda \beta \right) ^{-\alpha }={E}%
(e^{-\lambda X_{\infty }}\mid X_{0}=0)=\Phi (\lambda )\text{ as }%
t\rightarrow \infty .
\end{align*}
The function $\Phi _{t}(\lambda )$ is the PLSt of some rv $X_{t} $ with an
atom at $0$ with probability mass $e^{-\alpha t}$.

\section{Scaling}

\label{sec:scaling} In each example above, scaling the log-Laplace transform
according to 
\begin{equation}  \label{eq:rescale}
L(\lambda )\rightarrow L_{1}(\lambda )=\frac{1}{\sigma _{1}^{2}}L(\sigma
_{1}^{2}\lambda )
\end{equation}
defines the law of a new rv. Indeed, the scaled function $L_{1}(\lambda )$
defined from the LLt $L(\lambda )$ is itself the LLt of some rv if and only $%
L(\lambda )$ is the LLt of some infinitely-divisible rv (which is true in
our examples here).

Under this scaling \eqref{eq:rescale}, the mean $\mu =L_{1}^{\prime }(0)$
remains invariant while the variance is rescaled from $-L^{\prime \prime
}(0)=\sigma ^{2}$ to $-L_{1}^{\prime \prime }(0)=\sigma _{1}^{2}\sigma ^{2}$.
(For instance, a $\Gamma(\alpha,\beta)$ rv transforms to a 
$\Gamma(\alpha /A,\beta A)$ rv with $A=\sigma _{1}^{2}$.)
Then TL transforms according to 
\begin{equation*}
\sigma ^{2}=\mu ^{b}\rightarrow \sigma ^{2}=\sigma _{1}^{2}\mu ^{b}
\end{equation*}
where $\sigma _{1}^{2}$ (previously named $A$) is the variance of the new
scaled rv when its mean is $1$. The log-linear version of TL now includes a
non-zero constant term $a=\log \sigma _{1}^{2}$: 
\begin{equation*}
\log \sigma ^{2}=b\log \mu \rightarrow \log \sigma ^{2}=a+b\log \mu .
\end{equation*}

If $X$ has support $\mathbb{N}_{0}:=\{0,1,2,\ldots \}$, then under the scaling 
\eqref{eq:rescale} its support becomes $e^{a}\mathbb{N}_{0}:=\{0,e^{a},2e^{a},\ldots \}$. 
Scaling in the discrete case changes the original
'counting' support $\mathbb{N}_{0}$ to $e^{a}\mathbb{N}_{0}$,
emphasizing that $A=e^a$ plays the role of a
scale parameter of the underlying distribution. In this discrete dispersion
model, $A=e^a$ could be interpreted as the
individual mass of each individual constituting the population \cite{JK}.

If the support of $X$ is the whole real line or the nonnegative half line, the
above rescaling by $e^{a}$ does not change the support.

\section{Conclusion and open questions}

Empirical observations that sample means and sample variances are
approximately consistent with TL cannot specify the underlying distribution,
although they can reject some possibilities, as when the best estimates of $%
a\neq 0$ or $b\neq 1$ reject a Poisson distribution. As \cite{CohenTE}
observed, {for every $\mu >0$ and any rv $Y$ with mean $0$ and variance $1$,}
the rv $X=\mu +\mu ^{b/2}Y$ obeys TL with exponent $b$ and intercept $a=0$. $%
X$ can then be scaled following \eqref{eq:rescale} to include the affine
term $a$, resulting in a model obeying TL with any $a$ and $b$.

For each of the infinitely divisible two-parameter models in this note, we
showed that the TL exponent $b$ depends on the parameters in very specific
ways. We showed that the admissible range of $b$ depends on the
distribution, and in some cases on the parameters of the distribution, in
ways that could help identify the underlying law from data approximately
obeying TL, or at least could exclude some possible underlying laws.

The expression $b=(2-\alpha)/(1-\alpha)$ for the TL exponent in \eqref{bofalpha} arose in at least
three different prior examples of TL. First, TL holds asymptotically (as
sample size increases toward infinity) with exponent $b=(2-\alpha)/(1-\alpha)$ for
the sample variance and sample mean of samples from nonnegative stable
distributions with tail index $\alpha \in (0,1)$ 
\cite[page 663, Proposition 2]{BCD2017}. Such distributions have infinite
mean. Second, more generally, TL holds asymptotically with exponent $b=(2-\alpha)/(1-\alpha)$ for the sample variance and sample mean of samples of increasing
size from nonnegative distributions with regularly varying upper tails and
tail index $\alpha \in (0,1)$ \cite[page 6, their eq. (3.2)]{CDS2020}.
Again, such distributions have infinite mean. Third, TL holds asymptotically
with exponent $b=(2-\alpha)/(1-\alpha)$ for samples of increasing size from
nonnegative distributions with regularly varying upper tails and tail index $%
\alpha \in (0,1)$ when the sample variance is replaced by the sample upper
semivariance \cite[page 4, Theorem 2, their eq. 19]{BCTY2021}.

By contrast, for samples of increasing size from nonnegative distributions
with regularly varying upper tails and tail index $\alpha \in (0,1)$, TL
holds asymptotically, but with different formulas for the exponent $b$, when
the sample variance is replaced by the sample lower semivariance ($b=2$,
regardless of $\alpha \in (0,1)$), the sample lower local semivariance ($b=2$%
, regardless of $\alpha \in (0,1)$), or the sample local upper semivariance (%
$b=(2-\alpha ^{2})/(1-\alpha )$) \cite{BCTY2021}.

Can four appearances of (\ref{bofalpha}), one here and three earlier, be
coincidences? Or is some underlying process or mechanism common to all these
different appearances?

\section{Appendix}

Here we briefly sketch that if $X$ is SD, it is a weak limit of a pure-death
branching process with immigration (if $\mathbb{N}_{0}$-valued), or a weak
limit of a continuous-time L\'{e}vy-driven Ornstein-Uhlenbeck process (if
continuous). The analysis concerns the rvs $X$ obeying TL with $a=0$. From
the scaling transform \eqref{eq:rescale} introducing $a=\log \sigma _{1}^{2}$%
, the modifications for $a\neq 0$ could be readily obtained. In both
discrete and continuous cases, a population that is randomly annihilated is
randomly regenerated by the recurrent arrivals of random quantities of
immigrants, yielding a stationary or invariant distribution of population
size.

\subsection{Discrete self-decomposable rvs and pure-death branching
processes with immigration in continuous time}

van Harn et al. \cite{vanH} construct a regenerative process in continuous
time that produces discrete SD distributions in the long run. Consider a
continuous-time homogeneous compound Poisson process $P_{r}(t)$, $t\geq 0$, $%
P_{r}(0)=0$, having rate $r>0$, with pgf 
\begin{equation}
{E}_{P_{r}(0)=0}(z^{P_{r}(t)})=\exp \{-rt(1-h(z))\},  \label{f41}
\end{equation}
where $h(z)$, with $h(0)=0$, is the pgf of the sizes of the clones or
immigrant clusters arriving at the jump times of $P_{r}(t)$. Let 
\begin{equation}
\varphi _{t}(z)=1-e^{-t}(1-z)  \label{f41b}
\end{equation}
be the pgf of a pure-death branching process started with one particle at $%
t=0$. (More general subcritical branching processes could be considered.)
This expression of $\varphi _{t}(z)$ is easily seen to solve $\overset{.}{
\varphi }_{t}(z)=f(\varphi _{t}(z))=1-\varphi _{t}(z)$, $\varphi _{0}(z)=z$,
as is usual for a pure-death continuous-time Bellman-Harris branching
process \cite{Harris} with affine branching mechanism $f(z)=r_{d}(1-z)$ with
fixed death rate $r_{d}=1$. The distribution function of the lifetime of the
initial particle is $1-e^{-t}$. Let $X_{t}$ with $X_{0}=0$ be a random
process counting the current size of some population for which a random
number of individuals (determined by $h(z)$) immigrate at the jump times of $P_{r}(t)$. 
Each newly arrived individual is independently and immediately
subject to the pure death process above. We have 
\begin{equation}
\phi _{t}(z):={E}(z^{X_{t}})=\exp \left\{ -r\int_{0}^{t}[1-h(\varphi
_{s}(z))]ds\right\} ,\quad \phi _{0}(z)=1,  \label{f42}
\end{equation}
with $\phi _{t}(0)={P}(X_{t}=0)=\exp \{-r\int_{0}^{t}(1-h(1-e^{-s}))ds\},$
the probability that the population is extinct at $t$. As $t\rightarrow
\infty $, 
\begin{align}
\phi _{t}(z)\rightarrow \phi _{\infty }(z)& =\exp \left\{ -r\int_{0}^{\infty
}[1-h(1-e^{-s}(1-z))]ds\right\}  \notag  \label{f43} \\
& =\exp \left\{ -r\int_{z}^{1}\frac{1-h\left( u\right) }{1-u}du\right\} .
\end{align}
So $X:=X_{\infty },$ as the limiting population size of this pure-death
process with immigration, is a SD rv \cite{vanH}. Define the rv $X_{c}$
implicitly by requiring that $X\overset{d}{=}cX^{\prime }+X_{c}$, where $%
X^{\prime }$ is an iid copy of $X$ and $0<c<1$. Then 
\begin{equation*}
\phi _{X_{c}}(z)=\frac{\phi _{\infty }(z)}{\phi _{\infty }(1-c(1-z))}=\exp
\left\{ -r\int_{z}^{1-c(1-z)}\frac{1-h(u)}{1-u}du\right\}
\end{equation*}
is a pgf. In such models typically, a decaying subcritical branching
population is regenerated by a random number of incoming immigrants at
random Poissonian times.

\subsection{Continuous self-decomposable rvs and Ornstein-Uhlenbeck process
in continuous-time}

When $X$ is continuous and SD, $X$ is the limiting distribution of
population size as $t\rightarrow \infty$ of some Ornstein-Uhlenbeck  process $X_{t}$: 
\begin{equation*}
dX_{t}=-X_{t}dt+d\mathcal{L}_{t}, \ \ X_{0}=0,
\end{equation*}
driven by the L\'{e}vy process $\mathcal{L}_{t}$ for which 
\begin{equation*}
Ee^{-\lambda \mathcal{L}_{t}}=\Phi _{0}( \lambda ) ^{t}, \ \ t\geq 0,
\end{equation*}
where $\Phi _{0}( \lambda ) $ is the PLSt of an infinitely divisible rv
appearing in the representation (\ref{Phi}) of $\Phi ( \lambda )
=Ee^{-\lambda X}$. See \cite{Ju}.

We now show that for a TweBLE rv with $\alpha \in ( -\infty ,0) $, there is
no $L_{0}( \lambda ) =-\log \Phi _{0}( \lambda ) $ such that 
$L_{0}^{\prime }( \lambda ) $ is completely monotone on $( 0,\infty ) $.
This result means that the TweBLE rv for $\alpha \in ( -\infty ,0) $ is not
SD, just infinitely divisible. Indeed, with $L(\lambda )=-\log \Phi (\lambda)$, 
\begin{eqnarray*}
L_{0}(\lambda ) &=&\lambda L^{\prime }(\lambda )=(1-\alpha )^{1-\alpha
}\lambda (\theta +\lambda )^{\alpha -1}, \\
L_{0}^{\prime }(\lambda ) &=&(1-\alpha )^{1-\alpha }(\theta +\lambda
)^{\alpha -2}[\theta +\lambda \alpha ],
\end{eqnarray*}
with $L_{0}^{\prime }(\lambda )>0$\ only if $\lambda >\lambda _{c}=-\theta
/\alpha >0,$\ so not in the full range $\lambda \in (0,\infty )$. 
So $\Phi_{0}(\lambda )$ is not completely monotone on $(0,\infty )$, 
and is therefore not an infinitely divisible PLSt. 
This Poisson-gamma regime for which the limiting
distribution is a Poisson sum $P$ of iid gamma-distributed clusters of size 
$\Delta$ was studied by \cite[p. 17, section 3.3.2]{Eisler}, who underline 
what they call its ``impact inhomogeneity'': 
$E(\Delta )=C(\alpha )\cdot E(P)^{-1/\alpha }$ for
some constant $C(\alpha )>0$.

By contrast, if $\alpha \in (0,1)$, then $L_{0}(\lambda )$ may be written as 
\begin{equation*}
L_{0}(\lambda )=\int_{0}^{\infty }(1-e^{-\lambda x})\pi (x)dx
\end{equation*}
where 
\begin{equation*}
\pi (x)dx=\frac{1}{\Gamma (1-\alpha )}x^{-(\alpha +1)}(\alpha +\theta
x)e^{-\theta x}dx
\end{equation*}
is a tempered L\'{e}vy measure integrating $1\wedge x$. The driving process 
$\mathcal{L}_{t}$ of $X_{t}$, with PLSt $Ee^{-\lambda \mathcal{L}%
_{t}}=e^{-tL_{0}(\lambda )}$, is a subordinator and $X=X_{\infty }$ is a SD
TweBLE rv obtained as the limiting distribution of the corresponding Ornstein-Uhlenbeck process.
One could extend this construction to other SD subordinated L\'{e}vy families,
such as those in \cite{KIS, DA, JBe, BN, Schoutens}.

\section*{Acknowledgments}
We thank two reviewers for excellent constructive comments.
T.H. acknowledges partial support from the ``Chaire \textit{Mod\'{e}lisation
math\'{e}matique et biodiversit\'{e}'' }of Veolia-Ecole
Polytechnique-MNHN-Fondation X, and support from the labex MME-DII Center of
Excellence (\textit{Mod\`{e}les math\'{e}matiques et \'{e}conomiques de la
dynamique, de l'incertitude et des interactions}, ANR-11-LABX-0023-01
project). This work was also funded by CY Initiative of Excellence (grant ``%
\textit{Investissements d'Avenir}''ANR- 16-IDEX-0008), Project ``EcoDep''
PSI-AAP2020-0000000013. The ECODEP project is organized by Paul Doukhan. We
are grateful for this opportunity to collaborate.

\end{document}